   \documentclass[11pt]{article}

\begin{document}
   \newcommand{\be}{\begin{equation}}
   \newcommand{\ee}{\end{equation}}
   \newcommand{\tri}{\triangle}
   \newcommand{\e}{\epsilon}
   \renewcommand{\a}{\alpha}
   \renewcommand{\b}{\beta}
   \renewcommand{\k}{\kappa}
   \newcommand{\p}{\partial}
   \renewcommand{\phi}{\varphi}
   \newcommand{\pf}{{\noindent \sc Proof. }}
   \def\ca{{\cal A}}
   \newcommand{\ep}{\Box}
   \newcommand{\modulo}[1]{|#1|^2}
   \newcommand{\sob} [2]{H^{#2}_{#1}(0,1)}
   \newcommand{\Integers}{{\sf Z}\hspace{-1.7mm}{\sf Z}}
   \newcommand{\Naturals}{\rule{.15mm}{2.9mm}{\sf N}}
   \newcommand{\Reals}{\rule{.15mm}{2.4mm}{\sf R}}
   \newcommand{\integral}[1]{\int_{0}^{1}{#1}\;dx}
   \newcommand{\integralT}[1]{\int_{0}^{t}{#1}\;d\tau}
   \newcommand{\integers}{{\sf Z}\hspace{-1.4mm}{\sf Z}}
   \newcommand{\Complex}{\; \rule[.5mm]{.15mm}{2mm}\hspace{-1.4mm}{\sf C}}
   \newtheorem{lemma}{Lemma}[section]
   \newtheorem{proposition}[lemma]{Proposition}
   \newtheorem{theorem}[lemma]{Theorem}
   \newtheorem{corollary}[lemma]{Corollary}
   \newtheorem{remark}[lemma]{Remark}
   \newtheorem{definition}[lemma]{Definition}

 \renewcommand{\theequation}{\thesection.\arabic{equation}}
\catcode`\@=11
\@addtoreset{equation}{section}

\def\qed{\hbox{\vrule width 6pt height 6pt depth 0pt}}


   \title{Averaging Principle for   Quasi-Geostrophic  Motions\\
   under Rapidly Oscillating Forcing\thanks{This work is supported
   by a grant of NNSF of China 10001018 and NSF DMS-9973204.}}
   \author{Hongjun Gao$^1$    and \ \ Jinqiao Duan$^2$\\
\\
1.  Department of
   Mathematics\\Nanjing Normal University\\
   Nanjing 210097,  China\\
2. Department of Applied
   Mathematics\\Illinois Institute of Technology\\Chicago, IL 60616, USA\\
   E-mail: duan@iit.edu }
   \date{October 15, 2000}
   \maketitle
   \begin{abstract}
   In this paper, the averaging principle for   quasi-geostrophic
   motions with rapidly oscillating forcing   is
   proved, both on finite but large time intervals and on
  the entire time axis. This includes comparison estimate,
  stability estimate,  and convergence
  result between  quasi-geostrophic motions and its averaged motions.
  Furthermore,  the existence of almost periodic quasi-geostrophic
  motions and  attractor convergence   are  also investigated.

   \bigskip

   {\bf Key words}:  Quasi-geostrophic fluid flows, almost periodic
motions,  rapidly oscillating forcing,  averaging principle,
stable manifolds and unstable manifolds

\bigskip

   {\bf   Mathematics Subject Classification (2000):}
 34G20, 35Q35, 86A05, 76U05

\bigskip

  {\bf Abbreviated title:}  Averaging Principle for Quasi-Geostrophic Motions

\end{abstract}

\newpage


   \section{Introduction}

The   quasi-geostrophic(QG)   equation models
large scale geophysical flows. It is derived as an
approximation of the rotating Navier-Stokes  equations
by an asymptotic expansion in a small Rossby number. The
barotropic QG equation is written in terms of stream
function $\psi(x, y, t)$  (\cite{pedlosky}, 16, p. 234):
   \be
   \triangle\psi_t + J(\psi, \tri\psi) + \beta\psi_x = \nu\tri^2\psi -
   r\tri\psi + f(x, y, t),
\label{eqnpsi}
   \ee
   where   $\beta > 0$ the meridional
   gradient of the Coriolis parameter, $\nu > 0$ the viscous dissipation
   constant, $r > 0$ the Ekman dissipation constant,   $f(x, y, t)$ the
   wind forcing, and  $J(f, g) = f_x g_y - f_y g_x$ is the Jacobian operator.

     Equation (\ref{eqnpsi})   can be rewritten in terms of the relative vorticity
   $\omega(x, y, t) = \tri\psi(x, y, t)$ as
   \be
   \omega_t + J(\psi, \omega) + \beta\psi_x = \nu\tri\omega - r\omega +
   f(x, y, t).
\label{eqnomg}
   \ee
   For an arbitrary bounded planar domain $D$ with
sufficiently regular (such as, piecewise smooth)  boundary $\partial D $
   this equation can be supplemented by homogenous Dirichlet boundary
   conditions for both $\psi$ and $\omega$, namely, the no-penetration
and slip boundary conditions proposed by Pedlosky \cite{Ped96}, p. 34:
   \be
   \psi(x, y, t) = 0, \;\; \omega(x, y, t) = 0, \;\mbox{on}\; \partial D,
\label{bc}
   \ee
   together with an appropriate initial condition,
   \be
   \omega(x, y, 0) = \omega_0(x, y), \; \mbox{in}\; D.
\label{ic}
   \ee

    The global well-posedness (i.e., existence and uniqueness of smooth
   solution) of the dissipative model (\ref{eqnomg})-(\ref{ic})
  can be obtained similarly as
   in, for example, \cite{bct}, \cite{df}, \cite{wu}; see also \cite{bk}.
   Steady wind forcing has been used in numerical simulations(\cite{ci}).
   Brannan et al \cite{bdw} considered the effect of quasi-geostrophic
   dynamics under random forcing. Duan et al \cite{duan} and \cite{dk}
   obtained the existence of time periodic, time almost periodic
   quasi-geostrophic response of time periodic and time almost periodic
   wind forcing respectively.

In this paper, we assume that the right-hand
   side or the forcing term  of  the QG flow model
(\ref{eqnpsi})  is rapidly oscillating, i.e.,
it has the form
 $f(x,y, t) = f(x, y, \eta t) \stackrel{\textstyle \tri}{=}
   f(\eta t)$,  with parameter $\eta \gg 1.$
We also assume that $f$ has a time average.
With such forcing, it is desirable to understand the
fluid dynamics in some averaged sense,  and compare
the averaged flows with the original unaveraged flows.

The main result of this paper is the averaging principle for   
quasi-geostrophic motions with rapidly oscillating forcing,
     both on finite but large time intervals and on
  the entire time axis. This includes comparison estimate, 
  stability estimate,  and convergence
  result  (as  $\eta \to \infty$)  between  quasi-geostrophic 
  motions and its averaged motions.
 We also investigate the existence of almost periodic quasi-geostrophic
  motions under  almost periodic  forcing,  and
   the convergence of the attractor of the non-autonomous equation
 (\ref{eqnpsi}) to the attractor of the averaged
   autonomous equation as $\eta \to \infty$.

In \S 2, we study the averaging principle for the QG flow model
on finite  but large time intervals     and in \S 3, 
we extend the averaging principle to the  entire time axis.
In the rest of this section, we briefly review some background and provide
some preliminaries for later use.

Starting from the fundamental
   work of Bogolyubov \cite{bo} the averaging theory for ODE has been
   developed and generalized in a large number of works (see
   \cite{bm}--\cite{fi} and the references therein).
   Bogolyubov's main theorems have been generalized in \cite{dak}
   to the case of differential equations with bounded operator-valued
   coefficients. Some problems of averaging of differential equations with
   unbounded operator-valued coefficients have been considered in
   \cite{he}--\cite{sim} in the framework of abstract parabolic equations.
   In \cite{iiy}, Ilyin considered the averaging principle for
an equation of the form
   \be
   \partial_t u = N(u) + f(\eta t),
   \ee
   where $f$ is a given right-hand side and $\eta \gg 1$ is a large
   dimensionless parameter, and $f$ satisfies
   \be
   \lim\limits_{t\to \infty}\frac{1}{t}\int_0^tf(s) ds = f_0.
   \ee

   Standard abbreviations $L^2 = L^2(D)$, $H^k_0 = H_0^k(D), k =1, 2...,$
   are used for the common Sobolev spaces, with $<\cdot, \cdot>$ and
   $\|\cdot\|$ denoting the usual scalar product and norm, respectively, in
   $L^2$. We need the following properties and estimates  (see \cite{df})  of
   the Jacobian operator $J: H_0^1\times H_0^1\to L^1$:
   \be
   \int_DJ(f, g)h dx dy = - \int_DJ(f, h)g dx dy, \int_DJ(f, g)g dx dy = 0,
   \ee
   \be
   |\int_D J(f, gdx dy| \le \|\nabla f\|\|\nabla g\|,
   \ee
   for all $f, g, h \in H_0^1$, and
   \be
   |\int_DJ(\nabla f, g)\tri h dx dy| \le \sqrt{\frac{2|D|}{\pi}}\|\tri
   f\|\|\tri g\|\|\tri h\|,
   \ee
   for all $f, g, h \in H_0^2$.  We also recall the {\em Poincar\'{e}
   inequality} \cite{gt}
   \be
   \|g\|^2 \le \frac{|D|}{\pi}\int_D|\nabla g|dx dy =
   \frac{|D|}{\pi}\|\nabla g\|,
\label{poincare}
   \ee
   for $g\in H_0^1$, and  the {\em Young's inequality} \cite{gt}
   \be
   AB \le \frac{\e}{2}A^2 + \frac{1}{2\e}B^2,
   \ee
   where $A, B$ are nonnegative numbers and $\e > 0$.

   We can further rewrite  the QG flow model  (\ref{eqnomg}).   From
 \be
   \tri \phi = \omega, \;\;(x, y) \in D,\;\; \phi|_{\partial D} = 0,
\ee
   we   get  $\phi = \tri^{-1}\omega$.  Thus  (\ref{eqnomg})  can be
  rewritten as
   \be
   \omega_t + J(\tri^{-1}\omega, \omega) + \beta\partial_x\tri^{-1}\omega =
   \nu\tri\omega - r\omega + f(x, y, t).
   \ee
   Let
   $$- {\cal A} = \nu\tri - rI - \beta\partial_x\tri^{-1}.$$
  Then by a result in \cite{he}, we know that ${\cal A}$ is a sectorial
   operator, and hence  $e^{-{\cal A}t}$ generates an
analytic semigroup in $L^2$.

   We will give a sufficient condition to ensure the smallest eigenvalue of
   ${\cal A}$ to be positive. Consider the eigenvalue
equation $ {\cal A}u = \lambda u$.  We have the following energy
estimate
   $$\lambda \|u||^2 = \nu \|\nabla u\|^2 + r\|u\|^2 -
   \int_D\tri^{-1}u\partial_xu dxdy$$
   $$\ge \nu\|u\|^2 + r\|u\|^2 - \beta\|\tri^{-1}u\|\nabla u\| $$
  $$\ge
   \nu\|\nabla u\|^2 + r\|u\|^2 - \frac{\beta|D|}{\pi}\|u\|\|\nabla u\|
   $$
   $$\ge \nu\|\nabla u\|^2 + r\|u\|^2 - \frac{\beta|D|}{\pi}(a_1\|u\|^2 +
   a_2\|\nabla u\|^2)   $$
   \be  \ge (\nu - \frac{\beta|D|}{\pi}a_1)\|\nabla u\|^2 + (r -
   \frac{\beta|D|}{\pi}a_2)\|u\|^2,
   \ee
where we used the Poincar\'{e}
inequality  (\ref{poincare})  after the second inequality sign above,
and where arbitrary constants
$a_1, a_2$ satisfy $ a_1a_2 = \frac{1}{4}$.
 Therefore, when
   \be
   4\nu r > \frac{\beta^2|D|^2}{\pi^2},
\label{condition}
   \ee
 and if we take $a_2 = \frac{\beta|D|}{4r\pi}$,  then we have
   \be
   \lambda \|u||^2 \ge (\nu - \frac{\beta^2|D|^2}{4r\pi^2})\|\nabla u\|^2
   \ge (\nu - \frac{\beta^2|D|^2}{4r\pi^2})\frac{|D|}{\pi}\|u\|^2.
   \ee
   So, when $\nu, r, \beta$ and $|D|$ satisfy the condition
(\ref{condition}),  the
   smallest eigenvalue of ${\cal A}$ is positive. In this case,
the QG flow model
   $$
   \omega_t + {\cal A}\omega + J(\tri^{-1}\omega, \omega) = f(x, y, t)
   $$
   is a dissipative dynamical system.

  We note that the condition  (\ref{condition})   is
sharper  than  the corresponding condition in Duan and Kloeden \cite{dk}.
We  define the  fractional power of ${\cal A}$ as follows \cite{he}:
   $${\cal A}^{\alpha} = ({\cal A}^{-\alpha})^{-1}, \;\mbox{where}\;{\cal
   A}^{-\alpha} =
   \frac{1}{\Gamma(\alpha)}\int_0^{\infty}t^{\alpha-1}e^{-{\cal A} t}  dt.$$
   The corresponding domains $D({\cal A}^{\alpha})$ are Banach spaces with
   norm given by
   $$\|x\|_{\alpha}:= \|x\|_{D({\cal A}^{\alpha})} = \|{\cal
   A}^{\alpha}x\|.$$

For the rest of this section, we recall some definitions and
useful results for later use.

   \begin{theorem}\cite{he} The following estimates are valid:
   \be
   \|e^{-{\cal A}t}\|_{L^2\to L^2} \le Ke^{-at}, \;\;\;\;\;\;\;\;\; t \ge
   0,
   \ee
   \be
   \|{\cal A}^{\alpha}e^{-{\cal A}t}\|_{L^2\to L^2} \le
   \frac{K_{\alpha}}{t^{\alpha}}e^{-at}, \;\; t > 0,
   \ee
   where $K, K_{\alpha}$ are positive constants.
   \end{theorem}

   \begin{theorem}\cite{he}  Given two sectorial operators $A$ and $B$ in
   $L^2$, let $D(A) = D(B)$, $Re\sigma(A) > 0, Re\sigma(B) > 0$, and for
   some $\alpha\in[0, 1)$. Let the operator $(A-B)A^{-\alpha}$ be bounded
   in $L^2$. Then for every $\gamma\in[0, 1)$, $D(A^{\gamma}) =
   D(B^{\gamma})$, the two norms being equivalent.
   \end{theorem}

   \begin{definition}\cite{lez} A continuous function $f:\Reals\to X$ is
   called almost periodic(a.p.), if for every $\e > 0$ there exists a
   number $l = l(\e) > 0$ such that each interval $(T, T + l)$ contains a
   point $\tau = \tau_{\e}$(called an almost period) satisfying the
   inequality
   $$\mbox{Sup}_{t\in \Reals}\|f(t+\tau) - f(t)\| \le \e.$$
   \end{definition}
   If $f$ depends on other arguments, then the above inequality holds
   uniformly with respect to norms.

   It follows from the theory of a.p. functions that there exists a
   countable set of number $\lambda_{\alpha}$ for which
   $$\lim\limits_{T\to \infty}\frac{1}{2T}\int_{-T}^T
   f(t)e^{i\lambda_{\alpha}t} dt \not=0.$$
   The numbers $\{\lambda_{\alpha}\}$ are called the Fourier exponents of
   $f$\cite{lez}.

   \begin{definition}\cite{lez} A countable set of numbers
   $\{\omega_{\alpha}\}$ is called the frequency basis
(denoted by ${\cal M}_{f}$)  of an a. p. function
   $f$   if every $\lambda_{\alpha}$ can be
   uniquely respected by a linear combination of the numbers
   $\omega_{\alpha}$ with integer coefficients.
   \end{definition}

   \begin{definition}\cite{lez}  For a given  a.p. function $f$, the sequence
   $\{t_m\}$ is called $f-$current if
   $$\mbox{Sup}_{t\in \Reals} \|f(t + t_m) - f(t)\| \le \e_m\to
   0,\;\mbox{as}\; m\to \infty.$$
   \end{definition}

   \begin{theorem}\cite{lez} Given two almost periodic functions $f$ and
   $g$, suppose that every $f-$current sequence is also a $g-$current
   sequence. Then the frequency basis of $g$ is contained in that of $f$:
   ${\cal M}_{g}\subset{\cal M}_{f}$.
   \end{theorem}

We now turn to the averaging principle for the QG flow model.


   \section{Averaging Principle on Finite Time Intervals}

 In this section, we consider the  averaging principle for the QG flow model
on finite (but large) time intervals.   In the next section, we extend the result
to  the entire time axis.
We assume that the right-hand
   side or the forcing term  of  the QG flow model
(\ref{eqnpsi})  or  (\ref{eqnomg})   is rapidly oscillating, i.e.,
it has the form
 $f(x,y, t) = f(x, y, \eta t) =
   f(\eta t)$,  with parameter $\eta \gg 1.$
  Let $\eta \gg 1$ be a
   large dimensionless parameter. Setting
   $$ \tau = \eta t, \;\;\e = \eta^{-1},$$
   we obtain the equation in the so-called standard form
   \be
   \omega_{\tau} + \e{\cal A}\omega + \e J(\tri^{-1}\omega, \omega) = \e
   f(x, y, t).
   \ee
   We assume that $f$ has a time average in $D({\cal A}^{\gamma})$;
the value of
   $\gamma$ will be specified later on. More precisely, let $f(\tau), f_0
   \in {\cal A}^{\gamma}$ and suppose that
   \be
   \|{\cal A}^{\gamma}(\frac{1}{T}\int_t^{t+T}f(\tau) d\tau - f_0)\| \le
   \min(M_{\gamma}, \sigma_{\gamma}(T)),
   \ee
   where $M_{\gamma} > 0, \sigma_{\gamma}(T) \to 0, \;\mbox{as}\; T \to
   \infty$.

   We consider the averaged equation
   \be
   \bar\omega_{\tau} + \e{\cal A}\bar\omega + \e J(\tri^{-1}\bar\omega,
   \bar\omega) = \e f_0(x, y).
\label{eqnaverage}
   \ee
   By the method and result of \cite{bav}, \cite{hale} and \cite{te}, we
     obtain the semigroup $S_t$ corresponding to equation (\ref{eqnaverage})
    possesses absorbing sets in the space $H = L^2, V = D({\cal
   A}^{\frac{1}{2}}) = H_0^1$ and $D({\cal A})$, $\|\cdot\|$ and $\|\cdot\|_{D({\cal A})}$ denote the norm in $V$ and $D({\cal A})$. These sets are certain
   balls $B(R_0)$ in these spaces, where $R_0$ is large enough. This means
   that for every bounded set $B$
   $$S_tB \subset B(R_0), \;\mbox{for}\; t > t_0(B, R_0).$$
 In addition, the semigroup is uniformly bounded in these spaces, that
   is, given any ball, in particular, the ball $B(R_0)$, there exists a
   ball $B(R)$ such that
   $$S_tB(R_0) \subset B(R), \;\mbox{for}\; t > 0.$$
   By increasing $R$ we may assume that
   $$S_tB(R_0) \subset B(R- \rho), \;\mbox{for}\; t > 0, \rho > 0,$$
where $\rho$ is a positive constant.
   We first consider   the averaging principle in the space $V$. Given a
   point $\omega_0$ in $B_V(R_0)$,  we compare the  trajectories (solutions)
$ \omega(\tau) $, $\bar\omega(\tau)$ of
system (2.1) and (2.3)
   starting from this initial point.  Consider their difference on the interval
   $\tau\in [0, \frac{T}{\e}]$, $T$ being arbitrary but fixed. We suppose
   for the moment that $\omega(\tau) \in B_V(R)$. Then the difference $z(\tau) =
   \omega(\tau) - \bar\omega(\tau)$ satisfies the equation
   \be
   \partial_t z + \e {\cal A} + \e[J(\tri^{-1}\omega, \omega) -
   J(\tri^{-1}\bar\omega, \bar\omega)] = \e(f(\tau) - f_0(\tau)).
\label{eqnz}
   \ee

We first have some estimates on the nonlinear terms.
   \begin{lemma}
   The nonlinear operator $J(u, v)$ is a bounded Lipschitz map in the
   following sense:
   $$
   \|J(u_1, v_1) - J(u_2, v_2)\| \le
   $$
   \be C_{\frac{1}{2}}( \|u_1\|_{\frac{1}{2}} + \|u_2\|_{\frac{1}{2}} +
   \|v_1\|_{\frac{1}{2}} + \|v_2\|_{\frac{1}{2}})(\|u_1 -
   v_1\|_{\frac{1}{2}} + \|u_2 - v_2\|_{\frac{1}{2}}),
   \ee
   $$
   \|J(u_1, v_1) - J(u_2, v_2)\|_{\frac{1}{2}} \le
   $$
   \be
   C_0( \|u_1\|_{D(\ca)} + \|u_2\|_{D(\ca)} + \|v_1\|_{D(\ca)} +
   \|v_2\|_{D(\ca)})(\|u_1 - v_1\|_{D(\ca)} + \|u_2 - v_2\|_{D(\ca)}),
   \ee
where $C_{\frac{1}{2}}$ and $C_0$ are some positive constants.
   \end{lemma}
   \pf Since
   \begin{eqnarray}
   J(u_1, u_2) - J(v_1, v_2) &=& \nonumber\\
    (u_{1x} - v_{1x})u_{2y} + (u_{1y} - v_{1y})v_{1x} &-& (u_{2x} -
   v_{2x})u_{1y} - (u_{2y} - v_{2y})v_{2x},
   \end{eqnarray}
    (2.5) and (2.6) are obtained by direct estimtes.  Here the equivalence of norms $\|\cdot\|_{H^2}$ and $\|\cdot\|_{D(\ca)}$ is used.  $\hfill \qed$

Now we get back to equation (\ref{eqnz}).
   Inverting the linear operator we come to an equivalent integral equation
   \begin{eqnarray}
   z(\tau) &=& \e\int_0^{\tau}e^{-\e{\cal A}(\tau -s)}[J(\tri^{-1}\omega,
   \omega) - J(\tri^{-1}\bar\omega, \bar\omega)]ds \nonumber\\&+&
   \e\int_0^\tau e^{-\e{\cal A}(\tau -s)}(f(s) - f_0)ds.
   \end{eqnarray}
    Using (1.18) and (2.5), the $D(\ca^{\frac{1}{2}})$-norm of the first
   term in the right hand side satisfies the inequality
   $$\|\e\int_0^{\tau}\ca^{\frac{1}{2}}e^{-\e{\cal A}(\tau
   -s)}[J(\tri^{-1}\omega, \omega) - J(\tri^{-1}\bar\omega,
   \bar\omega)]ds\|
   $$
   $$\le \e\int_0^{\tau}K_{\frac12}\e^{-\frac12}(\tau -s)^{-\frac12}e^{-\e
   a(\tau - s)}2R\|z(s)\|_{\frac12}ds
   $$
   \be
   = 2RK_{\frac12}\e^{\frac12}\int_0^{\tau}(\tau -s)^{-\frac12}e^{-\e
   a(\tau - s)}\|z(s)\|_{\frac12}ds.
   \ee
   Let us estimate the second term in the right hand side of (2.8).
Integrating by parts we have
   $$
   \|\e\int_0^\tau e^{-\e{\cal A}(\tau -s)}(f(s) - f_0)ds\|_{\frac12}
   $$
   $$
   = \|-\e e^{-\e{\cal A}(\tau -s)}\int_s^{\tau}(f(t) - f_0)dt|_0^{\tau} +
   \e^2\int_0^\tau \ca e^{-\e{\cal A}(\tau -s)}\int_s^{\tau}(f(s) -
   f_0)ds\|_{\frac12}
   $$
   $$
   \le
   \|\e\ca^{\frac12-\gamma}e^{-\e{\cal
   A}\tau}{\ca}^{\gamma}\int_0^{\tau}(f(t) - f_0)dt\|
   $$
   \be
    + \|\e^2\int_0^\tau \ca^{\frac32-\gamma} e^{-\e{\cal A}(\tau
   -s)}\ca^{\gamma}\int_s^{\tau}(f(s) - f_0)ds\|.
   \ee
   Using (1.18) and (2.2), we  further have
   $$
   \|\e\ca^{\frac12-\gamma}e^{-\e{\cal
   A}\tau}{\ca}^{\gamma}\int_0^{\tau}(f(t) - f_0)dt\| \le
   \e K_{\frac12-\gamma}e^{-\e a\tau}(\e\tau)^{\gamma -
    \frac12}\|\frac{1}{\tau}\int_0^{\tau}{\ca}^{\gamma}(f(t) - f_0)dt\|
   $$
   $$
   = (\e\tau)^{\frac12 + \gamma}K_{\frac12-\gamma}e^{-\e
   a\tau}\|\frac{1}{\tau}\int_0^{\tau}{\ca}^{\gamma}(f(t) - f_0)dt\|
   $$
   \be
   \le (\e\tau)^{\frac12 + \gamma}K_{\frac12-\gamma}\min(M_{\gamma},
   \sigma_{\gamma}(\tau))e^{-\e a\tau} =: L(\tau).
   \ee
   For any $\delta > 0$, let $\tau_{\delta}$ be so large that for $\tau \ge
   \tau_{\delta}, \sigma_{\gamma} \le \delta$. Let $\e_0$ be so small that
   for $\e < \e_0$ then inequality $\frac{T}{\e} > \tau_{\delta}$ is valid.
   Then
   $$
   L(\tau) \le G_{\gamma 1}(T, \e) = e^{-\e a\tau}\left\{\begin{array}{cc}
   T^{\frac12 + \gamma}K_{\frac12 - \gamma}\delta, \;&\mbox{if} \;\tau \ge
   \tau_{\delta},\\
   (\e\tau)^{\frac12 + \gamma}K_{\frac12-\gamma}M_{\gamma}, \;&\mbox{if}\;
   \tau < \tau_{\delta}.
   \end{array}
   \right.
   $$
   Let $\gamma > - \frac12$. Since $\tau_{\delta}$ does not depend on $\e$,
   we let $\delta\to 0$ and then $\e\to 0$. We obtain
   \be
   \|\e e^{-\e{\cal A}\tau}\int_0^{\tau}(f(t) - f_0)dt\|_{\frac12} \le
   G_{\gamma 1}(T, \e) \to 0\;\mbox{when}\; \e \to 0.
   \ee
   $$
   \|\e^2\int_0^\tau \ca^{\frac32-\gamma} e^{-\e{\cal A}(\tau
   -s)}\ca^{\gamma}\int_s^{\tau}(f(s) - f_0)ds\| \le K_{\frac32 -
   \gamma}\e^{\frac12 + \gamma}\int_0^{\tau}\min(M_{\gamma},
   \sigma_{\gamma}(u))u^{\gamma - \frac12}du
   $$
   $$
   \le K_{\frac32 - \gamma}M_{\gamma}\e^{\frac12 +
   \gamma}\int_0^{\tau_{\mu}}u^{\gamma - \frac12}du + K_{\frac32 -
   \gamma}\e^{\frac12 + \gamma}\mu\int_0^{\frac{T}{\e}}u^{\gamma - \frac12}du
   $$
   \be
   = K_{\frac32 - \gamma}(\gamma + \frac12)^{-1}((\e\tau_{\mu})^{\frac12 +
   \gamma} + \mu T^{\frac12 + \gamma}) =: G_{\gamma 2}(T, \e),
   \ee
   where for any $\mu > 0$ we have chosen $\tau_{\mu}$ so large that
   $\sigma_{\gamma}(\tau) < \mu$ when $\tau > \tau_{\mu}$. Letting $\mu \to
   0$ and then $\e\to 0$ we obtain
   $$ G_{\gamma 2}(T, \e) \to 0, \;\; \e\to 0.$$
   Thus, by (2.8)--(2.13) we obtain the following inequality:
   \be
   \|z(\tau)\|_{\frac12} \le K\e^{\frac12}\int_0^{\tau}(\tau
   -s)^{-\frac12}\|z(s)\|_{\frac12}ds + G_{\gamma}(T, \e),
   \ee
   where $K = 2RK_{\frac12}$ and $G_{\gamma} = G_{\gamma 1} + G_{\gamma 2} \to
   0, \e \to 0$.

We need the following fact.
   \begin{lemma}\cite{he} Let $\gamma\in(0, 1]$ and for $t\in [0, T]$
   $$u(t) \le a + b\int_0^t(t-s)^{\gamma-1}u(s)ds.$$
   Then
   $$u(t) \le aE_{\gamma}((b\Gamma(\gamma))^{\frac{1}{\gamma}}t),$$
   where the function $E_{\gamma}(z)$ is monotone increasing and
   $E_{\gamma}(z)\sim \gamma^{-1}e^z$ as $z\to \infty$.
   \end{lemma}
    Applying this lemma to the inequality (2.14) on $\tau\in[0,
   \frac{T}{\e}]$, we obtain
   \be
   \|z(t)\|_{\frac12} \le G_{\gamma}(T, \e)E_{\frac12}(\e\tau\pi K^2) \le
   G_{\gamma}(T, \e)E_{\frac12}(T\pi K^2) = \eta_T^1(\e).
   \ee
   Using (1.18) and (2.6), we can do the same in $D(\ca)$ assuming
$\gamma > 0$ in (2.2) and obtain
   \be
   \|z(t)\|_{D(\ca)} \le F_{\gamma}(T, \e)E_{\frac12}(\e\tau\pi K^2)  =
   \eta_T^2(\e), \;F_{\gamma}(T, \e) \to 0\; \mbox{as}\; \e \to 0.
   \ee
 We thus have proved the proximity of solutions of (2.1) and (2.3)
 in $V$ and $D(\ca)$,  assuming that
   the  trajectory  $\omega(t)$    with  initial condition
$\omega(0) \in B_{V, D(\ca)}(R_0)$ stays in the ball $B(R)$ on the interval $[0, \frac{T}{\e}]$.

    Let $\e$ be so small that the right-hand side of (2.15) and (2.16) are
   less than $\frac{\rho}{2}$, where $\rho$ is defined earlier in this section
when we discuss absorbibg sets.
Suppose that the trajectory $\omega(t)$
   leaves the ball $B(R)$ during the interval $[0, \frac{T}{\e}]$ and let
   $\tau^*$ be the first momemt where $\|\omega(\tau^*) = R$. However, on
   the interval $\tau\in[0, \tau^*]$ both trajectories stay in the ball
   $B(R)$ and what we have proved so far shows that the inequality
   $\|\omega(\tau) - \bar\omega(\tau)\| \le \frac{\rho}{2}$ is valid. In
   particular, it is valid for $\tau = \tau^*$. This together with the
   inequality $\|\bar\omega(\tau^*)\| \le R - \rho$,  which holds by the
   hypothesis of the following theorem and the property of the semigroup
   $S(t)$, gives the contradiction
   $$\|\omega(\tau^*)\| \le \|\omega(\tau^*) - \bar\omega(\tau^*)\| +
   \|\bar\omega(\tau^*)\| \le R - \frac{\rho}{2}.
   $$

Therefore we have the following main result in this section.

   \begin{theorem}
(Averaging principle on finite time intevals)
   Let the right-hand side of equation (2.1) has an average in the sense of
   (2.2). Let $T > 0$ be arbitrary and fixed.

   If $\gamma > - \frac12$ and $\omega(0) = \bar\omega(0) \in B_V(R_0)$,
   that is, the initial values coincide and belong to the absorbing ball,
   then for $\tau\in[0, \frac{T}{\e}]$,
   $$
   \|\omega(\tau) - \bar\omega(\tau)\|_{\frac12} \le \eta_T^1(\e) \to 0\;
   \mbox{as}\; \e \to 0.
   $$
   If $\gamma > 0$, and $\omega(0) = \bar\omega(0) \in B_{D(\ca)}(R_0)$,
   then for $\tau\in[0, \frac{T}{\e}]$,
   $$
   \|\omega(\tau) - \bar\omega(\tau)\|_{D(\ca)} \le \eta_T^2(\e)\to 0\;
   \mbox{as}\; \e \to 0,
   $$
   where $\eta_T^1(\e)$ and $\eta_T^2(\e)$ are defined in (2.15), (2.16),
   respectively.
   \end{theorem}

This theorem gives comparison estimate and convergence result
(as $\eta  \to 0$)
between the QG flows and   averaged QG flows, on finite
but large time intevals.


   \section{Averaging Principle on the Entire Time Axis}

Now we turn to averaging principle for the QG flows on the entire time
axis.
   Consider
   \be
   \omega_{\tau} + \e{\cal A}\omega + \e J(\tri^{-1}\omega, \omega) = \e
   f(x, y, t).
   \ee
   All the hypotheses concerning the data of the problem are the same as
   those in \S 2; in particular, the average $f_0(x,y)$ 
   exists in the sense of (2.2).

We first consider  stationary {\em averaged} QG flows:
   \be
   \ca \omega_0 + J(\tri^{-1}\omega_0, \omega_0) = f_0.
   \label{stationary}
   \ee
Note that this is {\em not} the  stationary QG   flow model,
but the time-independent version of the {\em averaged} QG flow model.

    Under the condition of (1.15) and using Leray-Schauder fixed point
   theorem and elliptic regularity(\cite{gt}), we know the stationary 
 {\em averaged} QG flow model (\ref{stationary})  
   has a unique stationary solution. Note that this unique stationary solution
is denoted as $\omega_0 (x, y)$ and it is not to be confused with
 an initial  datum for the time-dependent QG flow model.

   We change the dependent variable (from $\omega$ to $z$)
   in equation (3.1) via:
   $$\omega = \omega_0 + z.$$
   Then by (3.2), we find that $z$ satisfies the equation
   $$\partial_t z = \e(- \ca z - J(\tri^{-1}\omega, \omega) +
   J(\tri^{-1}\omega_0, \omega_0) + f(\tau) - f_0).
   $$
   Since
   $$J(\tri^{-1}\omega, \omega) - J(\tri^{-1}\omega_0, \omega_0) =
   J(\tri^{-1}\omega, z) + J(\tri^{-1}z, \omega_0),
   $$
    we have
   $$
   \partial_t z = \e(- \ca z - J(\tri^{-1}\omega, z) - J(\tri^{-1}z,
   \omega_0) + f(\tau) - f_0).
   $$
   Changing the dependent variable once again  (from $z$ to $h$) via,
   $$z = h - \e v(\tau, \e),
   $$
   we obtain
   $$
   \partial_t z = \partial_t h - \e \partial_t v = \e(- \ca h + \e\ca v -
   J(\tri^{-1}\omega, h)
    - J(\tri^{-1}h, \omega_0)$$
   $$
   + (\e J(\tri^{-1}\omega, v) + \e J(\tri^{-1}v, \omega_0) + f - f_0).
   $$
   We chose  the auxiliary function $v(\tau, \e)$ to satisfy the equation
   \be
   \partial_{\tau}v = - \e{\ca}v + f_0 - f.
   \ee
   Then we obtain the following equation for the new dependent
variable $h$:
   $$
   \partial_{\tau}h = - \e({\ca}h  - J(\tri^{-1}\omega_0, h)
    - J(\tri^{-1}h, \omega_0)) +
   $$
   \be
    \e(- J(\tri^{-1}h, h) + J(\e v, h) + \e J(\tri^{-1}(\omega_0 - \e v),
   v) + \e J(\tri^{-1}h, v) + \e J(\tri^{-1}v , \omega_0)).
\label{eqnh}
   \ee

For the rest of the section, we study the equation (\ref{eqnh}) for
the new dependent variable $h$.
   Now,  we first consider equation (3.3) for the auxiliary
function $v(\tau, \e)$.

   \begin{lemma}
   Assume that the function $f$ has an average in the sense of (2.2).
If $\alpha -  \gamma < 1$, then equation (3.3) has a unique
solution $v(\tau, \e)$
   bounded in $D(\ca^{\alpha})$ uniformly in $\tau\in R$. Moreover,
   \be
   \|\e v(\tau, \e)\|_{\alpha} \to 0, \; \mbox{as}\; \e \to 0.
   \ee
   If $f$ is almost periodic with values in $D(\ca^{\gamma})$, then $v$ is
   almost periodic in $D(\ca^{\alpha})$ with frequency basis contained in
   that of $f$.
   \end{lemma}
   \pf The desired solution is given by the formula
   \be
   v(\tau, \e) = \int_{-\infty}^{\tau} e^{-\e\ca (\tau -s)}(f_0 -f))ds.
   \ee
   The uniqueness will be proved in a more general context in Lemma 3.2.
below.  Now
   we prove (3.5). Integrating by parts and using (2.2) and (1.18), we
   obtain
   \begin{eqnarray*}
   \|\e v(\tau, \e)\|_{\alpha} &=& \|\e \int_{-\infty}^{\tau} e^{-\e\ca
   (\tau -s)}(f_0 -f))ds\|_{\alpha}\\
   & =& \|\int_{0}^{\infty} e^{-\e\ca s}(f_0 - f(\tau -s))ds\|_{\alpha}\\
   &=& \|\e^2 \int_{0}^{\infty} \ca e^{-\e\ca s}\int_0^s(f_0 - f(\tau
   -t))dtds\|_{\alpha}\\
   &=& \|\e^2 \int_{0}^{\infty} \ca^{1+\alpha - \gamma} e^{-\e\ca
   s}s(s^{-1}\ca^{\gamma}\int_0^s(f_0 - f(\tau -t))dt)ds\|\\
   &\le& \e^{1+\alpha - \gamma}K_{1 + \alpha - \gamma}\int_0^{\infty}s^{\gamma
   - \alpha}e^{-\e a s}\min(M_{\gamma}, \sigma(s))ds\\
   &\le&K_{1 + \alpha - \gamma}(\e^{1+\alpha - \gamma}M_{\gamma}(1+\alpha -
   \gamma)^{-1}s_0^{\gamma - \alpha + 1} + \delta a^{1+\alpha -
   \gamma}\Gamma(1+\alpha - \gamma)),
   \end{eqnarray*}
   where $\delta$ is small and $s_0 = s_0(\delta)$ is so large that
   $\sigma(s) \delta$ when $s > s_0$. Letting $\delta \to 0$ and then $\e
   \to 0$, we obtain (3.5).

   Finally, let us prove the last statement of   Lemma 3.1.  By
   Theorem 1.6,  it is sufficient to show that every $f-$recurrent sequence
   $\{\tau_m\}$ is also $v-$recurrent. By (3.6),
   $$
   v(\tau + \tau_m) - v(\tau) = \int_0^{\infty}e^{-\e\ca s}(f(\tau-s) -
   f(\tau + \tau_m - s))ds.
   $$
   Therefore
   $$
   \mbox{sup}_{\tau \in R}\|v(\tau + \tau_m) -  v(\tau)\|_{\alpha} \le
   \mbox{sup}_{\tau \in R}\|f(\tau + \tau_m) - f(\tau)\|_{\gamma}K_{\alpha -
   \gamma}\Gamma(\gamma - \alpha + 1)(\e a)^{\alpha - \gamma }.
   $$
Thus $\{\tau_m\}$ is indeed $v-$recurrent.
   The proof  of Lemma 3.1 is complete.   $\hfill \qed$

 Now we go back to study the equation (\ref{eqnh}) for
the new dependent variable $h$.  We   consider the operator
   \def\cl{{\cal L}}
   \be
   {\cal L}h = {\ca}h  -  J(\tri^{-1}h, \omega_0) - J(\tri^{-1}\omega_0,
   h).
   \ee
   The operator $\cl + \lambda I$ has a compact inverse in $H = L^2$ for
   $\lambda > \lambda_0 = \frac{\pi}{2\lambda_1^3|D|}\|f_0\|^2 $, where
   $\lambda_1 \stackrel{\textstyle \tri}{=} \nu -
   \frac{\gamma^2|D|^2}{4r\pi^2}$.
In fact, we consider the equation
   $${\ca}h  -  J(\tri^{-1}h, \omega_0) - J(\tri^{-1}\omega_0, h) + \lambda
   h = f.
   $$
   Multiplying by $h$ in $L^2$ and observing that $(J(\tri^{-1}\omega_0, h),
   h) = 0$, we obtain
   $$(\ca h, h) - (J(\tri^{-1}h, \omega_0), h) + \lambda \|h\|^2 = (f, h).
   $$
   Using $H^2 \hookrightarrow L^{\infty}$, we find that
   $$
   - (J(\tri^{-1}h, \omega_0), h) \le \|\nabla h\|\|\nabla
   \omega_0\|_{\frac12}\|h\| \le \frac{\lambda_1}{2}\|\nabla h\|^2 +
   \frac{1}{2\lambda_1}\|h\|^2\|\nabla\omega_0\|^2.
   $$
   By (1.16) and (3.2) we have
   $$(\ca h, h) \ge \lambda_1\|\nabla h\|^2,\;\;\|\nabla\omega_0\|^2 \le
   \frac{\pi}{\lambda_1^2|D|}\|f_0\|^2.
   $$
   For $\lambda \ge  \lambda_0$, $((\cl  + \lambda)h, h)$ is coercive and
   by the Lax-Milgram lemma the equation
$$((\cl  + \lambda)h, h) = f$$
 has a
   unique solution $h\in D(\ca^{\frac12})$.  Note that the embedding $D(\ca^{\frac12})
   \to L^2$ is compact. Hence $\cl$ has a compact resolvent.
    Using (1.9),
   we can estimate $ J(\tri^{-1}h, \omega_0)$ as follows.
   $$| J(\tri^{-1}h, \omega_0) + J(\tri^{-1}\omega_0, h)| \le
   \|h\|\|\nabla \omega_0\| + \|\omega_0\|\|\nabla h\|.$$

   The operator $\ca$ is sectorial,  and due to the
   above estimates ,  we see that the operator
   $\cl$ is also  sectorial ( \cite{he}).
Moreover, if $Re \lambda > \lambda_0$, then the
   operator $\cl_{\lambda} = \cl + \lambda I$ is invertible; hence $$
   Re\;\sigma(\cl_{\lambda_0}) > 0.
   $$
   We also noe that
   $$
   (\cl_{\lambda_0} - \ca)\ca^{-\alpha} = {\lambda_0}\ca^{-\alpha}-
   J(\ca^{-\alpha}\tri^{-1}\cdot, \omega_0) - J(\tri^{-1}\omega_0,
   \ca^{-\alpha}\cdot).$$
   If $\alpha > 0$, the operator  $\cl_{\lambda_0}$
is bounded from $L^2 \to L^2$ (using
   (1.9)). The fact that $D(\ca) = D(\cl_{\lambda_0})$ follows from the
   theorem on the regularity of the 2-order elliptic problem[(\cite{gt}].
   It now follows from Theorem 1.2 that $D(\ca^{\alpha}) =
   D(\cl_{\lambda_0}^{\alpha}), \alpha\in [0, 1]$, that is
   \be
   c_{1\alpha}\|\cl_{\lambda_0}^{\alpha}\| \le \|\ca^{\alpha}\| \le
   c_{2\alpha}||\cl_{\lambda_0}^{\alpha}\|, \;\forall h \in
   D(\ca^{\alpha}).
   \ee
   We have thus shown that $\cl$ has a discrete spectrum. Let the stationary
   solution $\omega_0$ be such that $Re(\cl) \not= 0$ (which depends on the
   choice of $f_0$). In other words, we suppose that
   $$\sigma(\cl) = \sigma_+(\cl)\cup\sigma_-(\cl),
   \;\sigma_+(\cl)\cap\sigma_-(\cl) = \emptyset\;\mbox{and}
   \;Re\sigma_+(\cl) < -a, Re\sigma_-(\cl) > a, a > 0.
   $$
   We observe that $\sigma_+(\cl)$ is a finite set of eigenvalues. The $+$
   sign indicates the unstable modes and $-$ sign indicates the stable
   modes.

   Let $\gamma_+$ be a coutour in the left half-plane enclosing
   $\sigma_+(\cl)$. We set
   $$
   P_+ = \frac{1}{2\pi i}\int_{\gamma_+}(\lambda I - \ca)^{-1} d\lambda, \;
   P_- = I - P_+.
   $$
   The operator $\cl $ can be decomposed as
   $$
   \cl = \cl_+ + \cl_-, \; \cl_+ = P_+\cl, \; \cl_- = P_-\cl,
   $$
   where $P_+\cl \subset \cl P_+,\; P_-\cl \subset \cl P_-.$

   Setting $H_+ = P_+H, \mbox{dim} H_+ = N < \infty$ and $ H_- = P_-H$, we
   see that
   $$\cl_+H_+ \subset H_+,\;\;\cl_-H_- \subset H_-.$$

   The operator $\cl_+\in {\cal L}(H_+)$ (bounded linear operator space
   on $H_+$) and $-\cl_-$ generates an analytic semigroup in $H_-$.
Note that $P_+$
   and $P_-$ commute with $\cl$ and the semigroup $T(t) $ (which
   is generated by $-\cl$),
 in the sense that $P_+\cl \subset \cl P_+, P_-\cl
   \subset \cl P_-, P_+T(t) = T(t)P_+$ and $P_T(t) = T(t)P_-$ for $t \ge
   0$(\cite{kato}).

   We consider the following equation for $t\in \Reals$:
   \be
   \partial_t h + \cl h = f(t).
   \ee
   \begin{lemma}\cite{iiy}
   Let $f(t)\in L_{\infty}(\Reals; D(\cl^{\gamma})), \gamma \ge 0$. If
   $\alpha - \gamma < 1$, then the equation (3.9) has a unique solution $h$
   bounded in $D(\cl^{\alpha})$:
   \be
   \|h\|_{C_b(\Reals; D(\cl^{\alpha}))} \le K(\alpha,
   \gamma)\|f\|_{L_{\infty}(\Reals; D(\cl^{\gamma}))}.
   \ee
   \end{lemma}

We continue to study   equation (\ref{eqnh}) .
   Let $F(h, \e, \tau) = - J(\tri^{-1}h, h) + J(\e\; v, h) + \e
   J(\tri^{-1}(\omega_0 - \e\; v), v) + \e J(\tri^{-1}h, v) + \e
   J(\tri^{-1}v , \omega_0)$, then for $h_i\in D(\cl^{\frac12})$ and
   assuming that $\|h_i\|_{\frac12} \le \rho(i = 1, 2)$, we have
   \be
   \|F(h_1, \e, \tau) - F(h_2, \e, \tau)\| \le N_{\frac12}(\e, \rho)\|h_1 -
   h_2\|_{\frac12},
   \ee
   \be
   \|F(0, \e, \tau) \|  \le N_{\frac12}(\e),
   \ee
   where $N_{\frac12}(\e, \rho), \;N_{\frac12}(\e) \to 0$ as $\e, \rho\to
   0$.

   \noindent
   The proof of (3.11) and (3.12) can be obtained using {\bf Lemma 2.1} and
   {\bf Lemma 3.1}.   Moreover,  if
 $h_i\in D(\cl)(i = 1, 2)$, using {\bf Lemma 2.1} and
   {\bf Lemma 3.1}, we can prove $F: D(\ca)\to D(\ca^{\frac12})$ is a
   bounded Lipschitz map (using {\bf Theorem 1.2}).

   Converting  equation (\ref{eqnh}) to the original time variable
$t = \e\tau$, we obtain
   \be
   \partial_t h + \cl h = Q(h, \e, t),
\label{eqnh2}
   \ee
   where $Q(h, \e, t) = F(h, \e, \frac{t}{\e})$. Obviously, $Q(h, \e, t)$
   satisfies (3.11)  and  (3.12).

 \begin{lemma}
   Assume that $Q$ satisfies (3.11) and (3.12).
Then if $\e < \e_0$ for some $\e_0>0$ small enough,  (3.13) has a
   unique bounded solution  $h^*$ with the following properties:

   1. $\|h^*\|_{C_b(\Reals; D(\ca^{\frac12})} \le \delta(\e) \to 0$, as $\e
   \to 0$.

   2. There exists  an initial manifold ${\cal M}_-$ in $D(\ca^{\frac12})$,
   {\em codim}${\cal M}_- = N$, such that if $h_0 \in {\cal M}_-$ and
   $\|h_0\|_{\frac12} \le \rho$ for $\rho$ small enough, then the solution
   $h(t)$ of equation (3.13) with $h(0) = h_0$ satisfies the estimate
   \be
   \|h(t) - h^*(t)\|_{\frac12} \le Ke^{- a|t|}\|h_0 - h^*(0)\|_{\frac12},
   \; a > 0,\;\mbox{as}\; t \to \infty.
   \ee
   In particular, if $N = 0$, then the solution $h^*$ is asymptotically
   stable. There exists an initial manifold ${\cal M}_+$ in
   $D(\ca^{\frac12})$, $\dim {\cal M}_+ = N$, such that (3.14) holds as
   $t\to - \infty$.

   3. If the function $Q(h, \e, \cdot): \Reals \to D(\ca^{\frac12})$ is
   almost periodic, then the solution $h^*$ is almost periodic with
   frequency basis contained in that of $Q$.
   \end{lemma}
   This lemma is a special case of \cite{iiy}. Its proof is based on the
   contraction map principle, stability argument, {\bf Lemma 3.2} and {\bf
   Theorem 1.6}.

We now have the following main result  in this section.

   \begin{theorem} (Averaging principle on the entire time axis)
 Assume that the forcing
   $f$ on the right-hand side of the quasi-geostrophic flow model
(2.1) has an average in the sense of (2.2)
   uniformly in $t\in \Reals$.  Assume also that the spectrum 
   of the linear operator $\cl$  in (3.7) does not
    intersect  with the imaginary axis. Then for $\e < \e_0$
   small enough:

   1. In a small neighbourhood of the stationary {\em averaged} 
   quasi-geostrophic flow $\omega_0$,   
   the full quasi-geostrophic flow model
   (2.1) has a unique solution $\omega^*(\tau)$, which is bounded on the
   entire time axis and satisfies:
   $$\|\omega*(\tau) - \omega_0\|_{\frac12} \le \delta(\e), \;\mbox{as}\;
   \e \to 0,$$
   where $\delta(\e)$ is in Lemma 3.3.

   2. In the ball $B_{V}(\rho)\in H_-$ with $\rho$ small enough 
   and $\omega(0)\in B_{V}(\rho)$, there exists a stable 
   manifold ${\cal M}_-$ in $D(\ca^{\frac12})$, {\em codim}${\cal M}_- = N$, 
   such that
   if the initial condition $\omega(0) \in {\cal M}_-$ and $t \to \infty$,
  then
   \be
   \|\omega(t) - \omega^*(t)\|_{\frac12} \le ce^{- a|t|}\|\omega(0) -
   \omega^*(0)\|_{\frac12}.
   \ee
  There also exists an unstable manifold ${\cal M}_+$ in
   $D(\ca^{\frac12})$, $\dim{\cal M}_+ = N$, such that if
 the initial condition $\omega(0) \in
   {\cal M}_+$, then the inequality (3.15) holds as $t \to - \infty$. In
   particular, if $N = 0$, then the  unsteady quasi-geostrophic flow 
   $\omega^*$ is asymptotically  stable.

   3. If the forcing $f $ is almost periodic in
   $D(\ca^{\gamma}), \gamma > - \frac12$, then $\omega^*(\tau)$ is almost
   periodic in $D(\ca^{\frac12})$ with frequency basis contained in that of
   $f$.
   \end{theorem}

 \pf These assertions of the theorem follow from the representation
   $$\omega = \omega_0 + h(\tau, \e) - \e v(\tau, \e)$$
   and {\bf Lemmas 3.1, 3.3}.   $\hfill \qed$

This theorem gives stability estimate of stationary {\em averaged}
QG flows; stability conclusion of unsteady QG flows
near the stationary {\em averaged} QG flow; 
comparison estimate between unsteady QG flows; and  the
existence of almost time periodic QG motions under almost time periodic 
wind forcing, on the entire time axis.

 Combining this theorem and earlier discussion  in (1.15) and (1.16),
we have
   \begin{corollary}
   Under the assumption of {\bf Theorem 3.4} and
   $$ 4\nu r > \frac{\beta^2|D|^2}{\pi^2},\;\;\|f_0\| <
   \sqrt{\frac{2|D|}{\pi}}{\lambda_1^2},$$
where $\lambda_1 = \nu - \frac{\beta^2|D|^2}{4r\pi^2} $ as defined before,
   then for every $R > 0$ such that the  stationary {\em averaged}
   QG flow $\omega_0$ has norm $\|\omega_0\|_{\frac12}\le R$,  
   the following inequality is valid:
   $$\|\omega(\tau) - \omega^*(\tau)\|_{\frac12} \le C(R)e^{-\e
   a\tau}\|\omega(0) - \omega^*(0)\|_{\frac12}, \tau \to \infty,$$
   where $\omega(\tau) = \omega(\tau, \e)$ is the solution of  the full
   quasi-geostrophic flow model
   (2.1) with initial condition $\omega(0) = \omega_0$ and $\e < \e_0(R)$.
   \end{corollary}

   Combining {\bf Theorem 3.4} and {\bf Corollary 3.5}, 
   we conclude that there exists  an
   a. p.  solution  for (2.1) when $f$ is a. p.,
   under the assumptions of {\bf Corollary 3.5}. Here we give the
   restriction for $\|f_0\|$, not for $\|f\|$.  This result is obtained
   under different conditions from
those of Duan and Kloeden \cite{dk}.

   \def\ch{{\cal H}}

   Let $f$ be  a.p. with values in $H$. Then by the results of
   \cite{chv1}--\cite{chv3}, also see \cite{dk}, there exists a uniform
   attractor $A_{\ch} = A_{\ch(f(\eta\cdot))}$ of the non-autonomous
dynamical system (2.1),
   with $\ch = \{ f^{\tau}, f^{\tau}(t) = f(t + \tau),
   \tau\in\Reals\}_{C_b(\Reals;H)}$, and this attractor approaches to the
attractor of the averaged dynamical system (2.3):

   \begin{theorem} (Attractor convergence)
   Suppose that $f$ is a.p. in $H$, then
   $$ dist_{D(\ca^{\frac12})}(A_{\ch(f(\eta\cdot))}, \bar A) \to 0,
   \;\;   as \; \eta \to \infty,$$
   where $A_{\ch}$ is the attractor for the quasi-geostrophic flow
model (2.1), and $\bar A$ is the attractor of the averaged quasi-geostrophic
 flow model (2.3).
   \end{theorem}
   
This theorem claims that the uniform  attractor of the full
quasi-geostrophic flow model approaches to the attractor
of the averaged quasi-geostrophic flow model, as the parameter in
the oscillating forcing $\eta \to \infty$.


\section{Summary}

In this paper, we have discussed
  averaging principle for   quasi-geostrophic motions
under rapidly oscillating forcing  (characterizied 
by a large dimensionless parameter $\eta$),
     both on finite but large time intervals and on
  the entire time axis.  We have derived comparison estimate, 
  stability estimate   
and proved convergence
  result  (as $\eta \to \infty$ )  between  
  quasi-geostrophic motions and its averaged motions.

 We also investigated the existence of almost periodic quasi-geostrophic
  motions under  almost periodic  forcing,  and
   the convergence of the attractor of the  quasi-geostrophic  flow model
   to the attractor of the averaged quasi-geostrophic  flow model
     as $\eta \to \infty$.

\newpage

   \end{document}